\documentclass[12pt,leqno, twoside]{article}

\usepackage{amsfonts}
\usepackage{amsmath}
\usepackage{amssymb}

\newcommand{\R}{\mathbb{R}}
\newcommand{\Rn}{\mathbb{R}^n}

\newcommand{\Z}{\mathbb{Z}}
\newcommand{\Q}{\mathbb{Q}}

\newcommand{\dist}{\operatorname{dist}}
\newcommand{\rank}{\operatorname{rank}}
\newcommand{\sgn}{\operatorname{sgn}}
\newcommand{\I}{\operatorname{I}}

\newcommand{\ar}{\longrightarrow}
\newcommand{\inv }{^{-1}}

\newtheorem{theorem} {Theorem}[section]
\newtheorem{theorem*}{Theorem}
\newtheorem{prop*} {Proposition}
\newtheorem{lemma*}{Lemma}
\newtheorem{lemma}[theorem]{Lemma}
\newtheorem{cor}[theorem]{Corollary}
\newtheorem{cor*}{Corollary}
\newtheorem{prop}[theorem] {Proposition}

\newtheorem{definition*}{Definition}

\numberwithin{equation}{section}

\title{Immersions of spheres and algebraically constructible functions\thanks{%
							Iwona~Karolkiewicz, Aleksandra~Nowel and Zbigniew~Szafraniec\\
							University of Gda\'{n}sk,
              Institute of Mathematics \\
              80-952 Gda\'{n}sk, Wita Stwosza 57\\
              Poland\\
              Tel.: +48-58-5232059\\
              Fax: +48-58-3414914\\
              Email: ikarolki@manta.univ.gda.pl \\
              Email: Aleksandra.Nowel@math.univ.gda.pl\\
              Email: Zbigniew.Szafraniec@math.univ.gda.pl\\}\thanks{%
							First and second named authors supported by the grant BW /5100-5-0286-7}
}

\author{Iwona~Karolkiewicz \and Aleksandra~Nowel \and Zbigniew~Szafraniec}

\date{April 2007}

\begin{document}

\def\nothanksmarks{\def\thanks##1{\protect\footnotetext[0]{\kern-\bibindent##1}}}
        
\nothanksmarks        

\maketitle

\pagestyle{myheadings}

\markboth{\textsc{\small I.~Karolkiewicz, A.~Nowel, Z.~Szafraniec}}{\textsc{\small Immersions of spheres and algebraically constructible functions}}

\begin{abstract}
Let $\Lambda$ be an algebraic set and let $g:\R ^{n+1}\times \Lambda \ar \R^{2n}$ ($n$ is even) be a polynomial mapping such that for each $\lambda\in \Lambda$ there is $r(\lambda)>0$ such that the mapping $g_{\lambda}=g(\cdot ,\lambda)$ restricted to the sphere $S^n(r)$ is an immersion for every $0<r<r(\lambda)$, so that the intersection number $\I (g_{\lambda}|S^n(r))$ is defined. Then $\Lambda \ni \lambda \mapsto \I (g_{\lambda}|S^n(r))\in \Z$ is an algebraically constructible function.
\end{abstract}

{\small
\noindent
{\bf Keywords:} immersions of spheres, algebraically constructible functions, real
algebraic sets.

\noindent
{\bf Mathematics Subject Classification (2000):} 14P25  57N35  32S50}

\section{Introduction} \label{sec:1}

McCrory and Parusi\'{n}ski \cite{parus-mccr1} have introduced algebraically constructible functions in order to study the topology of real algebraic sets. 

Let $\Lambda$ be a real algebraic set. An integer valued function $\phi:\Lambda \ar \Z$ is {\em algebraically constructible} if there exist an algebraic set $W$ and a proper regular map $p:W\ar \Lambda$ such that $\phi(\lambda)$ equals the Euler characteristic $\chi (p\inv (\lambda))$.

If that is the case then $\phi$ is {\em semialgebraically constructible}, that is there exists a semialgebraic stratification ${\cal S}$ of $\Lambda$ such that $\phi$ is constant on strata of ${\cal S}$. If $\Lambda$ is irreducible then $\phi$ has to be generically constant modulo $2$, see for instance \cite[Proposition 2.3.2]{ak-ki2}, and there exist a real polynomial $g:\Lambda \ar \R$ and a constant $\mu$ such that generically on $\Lambda$: $\phi \equiv \mu +\sgn g \mod 4$ \cite{costekurdyka2}.

In fact, algebraically constructible functions are precisely those constructible functions which are sums of signs of polynomials \cite{szafr-parus1}, \cite{szafr-parus2}.

Let $f:\Rn \times \Lambda\ar \Rn$ be a polynomial mapping such that $0\in \Rn$ is isolated in $f(\cdot,\lambda)\inv (0)$ for all $\lambda \in \Lambda$, so that the local topological degree $\deg _0f(\cdot, \lambda)$ of $f(\cdot, \lambda)$ at the origin is well defined. Then $\Lambda \ni \lambda\mapsto \deg _0f(\cdot, \lambda)\in \Z$ is algebraically constructible \cite[Theorem 3.4]{szafr-parus1}.

For a recent account of the theory we refer the reader to \cite{coste1}, \cite{parus-mccr2}, \cite{parus-mccr3}.

\medskip

Whitney \cite{whitney1} has introduced an intersection number $\I (g)$ for an immersion $g:M^n \ar \R^{2n}$. If $n$ is even then $\I(g)\in \Z$, if $n$ is odd then $\I(g)\in \Z _2$. Smale \cite{smale2} proved, that two immersions $f,g:S^n\ar \R^{2n}$ are regularly homotopic if and only if $\I (f)=\I (g)$. 

In this paper we investigate how does the ``local'' intersection number change when there is an algebraic family of immersions. Let $g=(g_1,\ldots ,g_{2n}):\R ^{n+1}\times \Lambda \ar \R ^{2n}$ ($n$ is even) be a polynomial mapping. Assume that for each $\lambda \in \Lambda$ there exists $r(\lambda)>0$ such that the mapping $g_{\lambda}=g(\cdot, \lambda)$ restricted to the sphere $S^n(r)$ is an immersion for $0<r<r(\lambda)$. If that is the case then the intersection number $\I (g_{\lambda}|S^n(r))$ is the same for all $0<r<r(\lambda)$.
We shall prove (Theorem~\ref{alg_constr}) that the function $$\Lambda \ni \lambda \mapsto \I (g_{\lambda}|S^n(r))\in \Z$$ is algebraically constructible.

\section{Preliminaries} \label{sec:2}

Suppose that $L$ is a $p$--dimensional oriented manifold, $H:L\ar \R ^p$ is a smooth mapping, and $U$ is an open subset of $L$ such that $H\inv (0)\cap U$ is compact.

There exists $(N, \partial N)$ --- a compact $p$--dimensional oriented manifold with boundary such that $N\subset U$ and $H\inv (0)\cap U\subset N\setminus \partial N$. 

By {\em the topological degree} of the mapping 
$$(U,U\setminus H\inv (0))\ni x \mapsto H(x)\in (\R ^p,\R ^p \setminus \{ 0\})$$
we mean the topological degree of
$$(N,\partial N)\ni x \mapsto H(x)\in (\R ^p,\R ^p \setminus \{ 0\})$$
which equals the degree of the mapping
$$\partial N \ni x \mapsto \frac{H(x)}{\| H(x)\|}\in S^{p-1}.$$
Of course, the degree does not depend on the choice of $N$.

\bigskip

Let $F_1,\ldots ,F_k ;G_1,\ldots ,G_{n-k}: \R^n,0 \rightarrow \R ,0$ be analytic functions 
defined in a neighbourhood of the origin.

Denote
\begin{align*}
F & =(F_1,\ldots,F_k):\R^n,0\rightarrow\R^k,0 \\
G & =(G_1,\ldots,G_{n-k}):\R^n,0\rightarrow\R^{n-k},0 \\
S^{n-1}(r) & =\{x\in\R^n\ |\ \| x\|=r\} \\
S^{n-1} & =\{x\in\R^n\ |\ \| x\|=1\} \\
B^n(r) & =\{x\in\R^n\ |\ \| x\|\leq r\}
\end{align*}

Suppose that $F^{-1}(0)$ has an isolated singularity at the origin,
i.e. for $\| x\|$ small enough, if $F(x)=0$ and $\rank [DF(x)]<k$
then $x=0$.
If $r>0$ is small enough then $S^{n-1}(r)$ cuts $F^{-1}(0)$ transversally,
so $M(r)=S^{n-1}(r)\cap F^{-1}(0 )$ is either void or a compact
$(n-k-1)$--dimensional manifold.

We shall say that vectors $v_1,\ldots,v_{n-k-1}$ in the tangent space
$T_x M(r)$ are {\em well oriented} if
$\nabla F_1(x),\ldots,\nabla F_k(x),x,v_1,\ldots,v_{n-k-1}$
are well oriented in $\R^n$. This way $M(r)$ is oriented.

Let $y\in \R^k$ be a regular value of $F$. Then $F^{-1}(y)$
is either void or an $(n-k)$--dimensional manifold. We shall say that
$w_1,\ldots,w_{n-k}\in T_x F^{-1}(y)$ are {\em well oriented}
if $\nabla F_1(x),\ldots,\nabla F_k(x),w_1,\ldots,w_{n-k}$ are well oriented in $\R^n$.
This way $F^{-1}(y)$ is also oriented.

Fix small $r>0$.
If $y$ lies sufficiently close to the origin then $S^{n-1}(r)$ cuts $F^{-1}(y)$
transversally, so $\overline{M}(r)=S^{n-1}(r)\cap F^{-1}(y)$ is either void or
a compact $(n-k-1)$--dimensional manifold. Moreover, $B^n(r)\cap F^{-1}(y)$
is a compact oriented manifold with boundary $\partial\, (B^n(r)\cap F^{-1}(y))=\overline{M}(r)$.
The orientation of the boundary induced from $B^n(r)\cap F^{-1}(y)$ may be described
the same way as the orientation of $M(r)$. Manifolds $M(r)$ and $\overline{M}(r)$ are isotopic.

Suppose that $F^{-1}(0)\cap G^{-1}(0 )=\{0\}$. Then
$M(r)\cap G^{-1}(0)=\emptyset$ and
$\overline{M}(r)\cap G^{-1}(0)=\emptyset$.

If $\| y\|$ is small enough then mappings
$$M(x)\ni x \mapsto \frac{G(x)}{\| G(x)\|}\in S^{n-k-1},$$
$$\overline{M}(x)\ni x \mapsto \frac{G(x)}{\| G(x)\|}\in S^{n-k-1}$$
have the same topological degree. Denote it by $\deg(M(r))$.
It equals the degree $\rho$ of $\overline{G}$, where $\overline{G}$ is the restricted mapping
$$(F^{-1}(y)\cap B^n(r),\overline{M}(r))\ni x\mapsto G(x)\in (\R^{n-k},\R^{n-k}\setminus\{0\}).$$
Choose a regular value $z\in\R^{n-k}$ near the origin. Then
$$\rho=\sum \sgn \det [D\, \overline{G}(x)]\ \ \ \ \ \ \mbox{(where }x\in \overline{G}^{-1}(z)\mbox{)}.$$
Using the concept of the Gram determinant the reader may check that
$$\rho=\sum \sgn \det [D(F,G)(x)],$$
where $x\in F^{-1}(y)\cap B^n(r)\cap G^{-1}(z)=(F,G)^{-1}(y,z)\cap B^n(r)$.
If $0<\| z\| \ll \| y\|\ll r\ll 1$ then $\rho$ equals the local topological degree $\deg_0(F,G)$
at the origin,
i.e. the topological degree of the mapping
$$S^{n-1}(r)\ni x\mapsto \frac{(F(x),G(x))}{\| (F(x),G(x))\|}\in S^{n-1}.$$

Set
$$N(r)=\{x\in S^{n-1}(r)\ |\ F_1(x)=\ldots=F_{k-1}(x)=0, F_k(x)\geq 0\}.$$
Suppose that $(F_1,\ldots,F_{k-1})^{-1}(0)$ has an isolated singularity at the origin.
Since $F^{-1}(0)$ has an isolated singularity, $N(r)$ is either void or
an $(n-k)$--dimensional compact oriented manifold with boundary, and
$\partial\, N(r)=S^{n-1}(r)\cap F^{-1}(0)=M(r)$. There is the restricted function
$$(N(r),\partial\, N(r))\ni x\mapsto G(x)\in (\R^{n-k},\R^{n-k}\setminus\{0\}),$$
which will be denoted by $\overline{G}_r$.

The topological degree of $\overline{G}_r$ equals the topological degree of
$$\partial\, N(r)\ni x\mapsto \frac{G(x)}{\|G(x)\|}\in S^{n-k-1},$$
which is the same as
$\deg(M(r))= \rho=\deg_0(F,G)$. We have proved:

\begin{prop} \label{loctopdeg}
Let $F_1,\ldots ,F_k ;G_1,\ldots ,G_{n-k}: \R^n,0 \rightarrow \R ,0$ be analytic functions
defined in a neighbourhood of the origin.
Suppose that $(F_1,\ldots,F_{k-1})^{-1}(0)$, as well as $(F_1,\ldots,F_{k})^{-1}(0)$, 
has an isolated singularity at the origin, and $\{0\}$ is isolated in
$(F_1,\ldots ,F_k ;G_1,\ldots ,G_{n-k})^{-1}(0)$.

If $r>0$  is small enough then the topological degree of the mapping $\overline{G}_r$, i.e. of
$$(N(r),\partial\, N(r))\ni x\mapsto G(x)\in (\R^{n-k},\R^{n-k}\setminus\{0\}),$$
does not depend on $r$ and is equal to the local topological degree
$\deg _0(F,G)$. 

\hfill $\Box$
\end{prop}

Let $\delta:\R^n,0\rightarrow \R,0$ be an analytic function. In \cite{szafraniec4}
there is proven that if $\alpha$ is a sufficiently large positive  even integer and $t\neq 0$ 
then $\delta-t\| x\|^\alpha$ has an isolated critical point at the origin. Using the same arguments
one may prove: 

\begin{prop} \label{isolsing}
Suppose that $(F_1,\ldots,F_{k-1})^{-1}(0)$ 
has an isolated singularity at the origin. 

Then there exists $\alpha_0>0$ such that for any even integer $\alpha > \alpha_0$ and $t\neq 0$
$$(F_1,\ldots,F_{k-1},\delta-t\| x\|^\alpha)^{-1}(0)$$
has an isolated singularity at the origin. 

\hfill $\Box$
\end{prop}

Suppose that $\delta\geq 0$ and 
$$X:=\{0\}\cup ((F_1,\ldots,F_{k-1},G_1,\ldots,G_{n-k})^{-1}(0)\setminus \delta^{-1}(0))$$
is closed. 
Then $\delta^{-1}(0)\cap X=\{0\}$.
Since $X$ is closed semianalytic, there exists a {\L}ojasiewicz exponent $\alpha_1$, such that
$\delta(x)>t\| x\|^{\alpha_1}$ for $t>0$ and $x\in X\setminus\{0\}$ sufficiently close to the origin.

Set
$$L(r)=\{x\in S^{n-1}(r)\ |\ F_1(x)=\ldots=F_{k-1}(x)=0\}.$$
Then $U(r)=L(r)\setminus\delta^{-1}(0)$ is an open subset of $L(r)$, and
$U(r)\cap G^{-1}(0)=S^{n-1}(r)\cap X$ is a compact subset of $U(r)$.

Let $\alpha > \max(\alpha_0,\alpha_1)$ be an even positive integer. Set $F_k=\delta-t\| x\|^\alpha$ for arbitrary $t>0$.
Then $F_k(x)>0$ for $x\in L(r)\cap G^{-1}(0)\setminus\delta^{-1}(0)$,
and $F_k(x)<0$ for $x\in L(r)\cap G^{-1}(0)\cap\delta^{-1}(0)$.

\begin{prop} \label{topdeg}
Suppose that \begin{itemize}
\item[(a)] $(F_1,\ldots,F_{k-1})^{-1}(0)$ has an isolated singularity at the origin,
so that $L(r)$ is either void or a compact oriented $(n-k)$--dimensional manifold
for small $r>0$,
\item[(b)] $\delta\geq 0$ and $U(r)\cap G^{-1}(0)=L(r)\cap G^{-1}(0)\setminus\delta^{-1}(0)$ is a compact
subset of $U(r)$.
\end{itemize}
Then there is a mapping $G_r$ given by
$$(U(r),U(r)\setminus G^{-1}(0))\ni x\mapsto G(x)\in (\R^{n-k},\R^{n-k}\setminus\{0\})$$
such that if $\alpha>\max(\alpha_0,\alpha_1)$ is  an even integer, $t>0$, and $F_k=\delta-t\| x\|^\alpha$, then
$$\R^n\ni x\mapsto (F(x),G(x))\in\R^n$$ has an isolated zero at the origin, and for each
$r>0$ small enough the topological degree of $G_r$ equals the local topological
degree $\deg_0(F,G)$.
\end{prop}
{\em Proof.} Since $F_k(x)>0$ for $x\in U(r)\cap G^{-1}(0)=L(r)\cap G^{-1}(0)\setminus\delta^{-1}(0)$,
then  $U(r)\cap G^{-1}(0)\subset N(r)\setminus \partial\, N(r)\subset U(r)$.
So the topological degree of $G_r$ equals the topological degree of $\overline{G}_r$, and by Proposition~\ref{loctopdeg}, equals $\deg_0(F,G)$. 

\hfill $\Box$

\section{The intersection number of an immersion} \label{sec:3}

Let $M$ be an $n$--dimensional manifold.
A $C^1$ map $g:M\ar \R ^m$ is called {\em an immersion} if for each $p\in M$ the rank of $Dg(p)$ equals $n$.

\medskip

A homotopy $h_t:M\ar \R ^m$ is called {\em a regular homotopy,} if at each stage it is an immersion and the induced homotopy of the tangent bundle is continuous.

\begin{theorem} \label{reg-homot}
{\em \cite[Theorem B]{smale2})}
Two $C^{\infty}$ immersions from an $n$--dimensional sphere $S^n$ to $\R ^m$ are regularly homotopic when $m\geqslant 2n+1$.
\end{theorem}

As in \cite{whitney1} we say that an immersion $g:M\ar \R^{2n}$ has a {\em regular self--intersection} at the point $g(p)=g(q)$ if 
$$Dg(p)T_pM+Dg(q)T_qM=\R ^{2n}.$$
An immersion $g:M\ar \R^{2n}$ is called {\em completely regular} if it has only regular self--intersections and no triple points.

\bigskip

Assume that $n$ is {\bf even} and $M$ is compact and oriented.

\medskip

Let $g:M\ar \R^{2n}$ be a completely regular immersion and have a regular self--intersection at the point $g(p)=g(q)$. 

Let $u_1,\ldots ,u_n\in T_pM$, $v_1,\ldots ,v_n\in T_qM$ be sets of well--oriented, independent vectors in respective tangent spaces of $M$. Then the vectors $Dg(p)u_1,\ldots , Dg(p)u_n$, $Dg(q)v_1,\ldots ,Dg(q)v_n$ form a basis in $\R ^{2n}$. As in \cite{whitney1} we will say that the self--intersection at the point $g(p)=g(q)$ is {\em positive} or {\em negative} according to whether this basis determines the positive or negative orientation of $\R^{2n}$.

{\em The intersection number} of a completely regular immersion $g$ is the algebraic number of its self--intersections.

\medskip

Let for a given (not necessarily completely regular) immersion $g:M\ar \R^{2n}$ define $G:M\times M \ar \R^{2n}$ as $G(x,y)=g(x)-g(y)$. Set 
$$\Delta =\{ (p,p) \ |\ p\in M\} \subset M\times M.$$
Since $g$ is an immersion, $\Delta$ is isolated in $G\inv (0)$, and so $G\inv (0)\setminus \Delta$ is a compact subset of $M\times M\setminus \Delta$. Of course, $M\times M\setminus \Delta$ is an open subset of $M\times M$, and $$(M\times M\setminus \Delta)\setminus G\inv (0)=(M\times M\setminus (G\inv (0)\cup \Delta)).$$ The topological degree $d(g)$ of the mapping
$$(M\times M\setminus \Delta, M\times M\setminus (G\inv (0)\cup \Delta))\ni (x,y)\mapsto G(x,y)\in (\R^{2n},\R^{2n}\setminus \{ 0\})$$  
is always an even integer. Let us denote
$$\I(g)=\frac{1}{2}d(g).$$
By  \cite[Theorem 3.1]{lashof-smale1} we have

\begin{theorem} \label{ind-deg} If $g$ is a completely regular immersion then its intersection number is equal to $\I(g)$.
\end{theorem}

\medskip

Assume that either $n$ is {\bf odd} or $M$ is non--orientable.

In this case one can also define {\em the intersection number} of a completely regular immersion $g:M\ar \R^{2n}$ as the number of its self--intersections modulo $2$. 

\medskip

By \cite[Theorem 2]{whitney1} if $M$ is closed then the intersection number is invariant under regular homotopies. 
As in \cite{whitney3}, if $M$ is closed, any immersion $g:M\ar \R^{2n}$ 
can be made completely regular by a regular homotopy.
So in this case we can define the intersection number also for $g$ that is not completely regular.

We have a characterisation of regularly homotopic immersions due to Smale \cite{smale2}:

\begin{theorem} \label{smale}
{\em \cite[Theorem C]{smale2}}
Two $C^{\infty}$ immersions $f$, $g$ from an $n$--dimensional sphere $S^n$ to $\R^{2n}$ are regularly homotopic if and only if $\I (f)=\I (g)$.
\end{theorem}

\section{Immersions on small spheres} \label{sec:4}

Let
$$h=(h_1,\ldots,h_l) :\R ^{n} \ar \R ^{l}$$ 
$$g=(g_1,\ldots,g_k) :\R^ {n} \ar \R ^{k}$$
be $C^1$ mappings. Put $M:=h^{-1}(0)$. Suppose that each point $p\in M$ is a regular point of $h$, i.e. the rank of the derivative matrix $Dh(p)$:
$$\left [
\begin{array}{cccc}
   \frac{\partial h_1}{\partial x_1}(p) & \frac{\partial h_1}{\partial x_2}(p)&\ldots & \frac{\partial h_1}{\partial x_n}(p)\\
   \vdots & \vdots & \ddots & \vdots\\
   \frac{\partial h_l}{\partial x_1}(p) & \frac{\partial h_l}{\partial x_2}(p)&\ldots & \frac{\partial h_l}{\partial x_n}(p)\\
 \end{array} \right ]$$ 
equals $l$ at each $p\in M$. If that is the case then $M$ is a $C^1$ $(n-l)$--manifold, and there is the restricted mapping $g|M:M\ar \R ^k$.

\begin{prop}\label{rank}
$$\rank
\left [
\begin{array}{c}
   D(g|M)(p)
\end{array} \right ]
=\rank
\left [
\begin{array}{ccc}
   \frac{\partial g_1}{\partial x_1}(p) &\ldots & \frac{\partial g_1}{\partial x_n}(p)\\
   \vdots & \ddots & \vdots\\
   \frac{\partial g_k}{\partial x_1}(p) &\ldots & \frac{\partial g_k}{\partial x_n}(p)\\
   \frac{\partial h_1}{\partial x_1}(p) &\ldots & \frac{\partial h_1}{\partial x_n}(p)\\
   \vdots & \ddots & \vdots\\
   \frac{\partial h_l}{\partial x_1}(p) &\ldots & \frac{\partial h_l}{\partial x_n}(p)\\
\end{array} \right ]
-l$$
at each point $p\in M$.
\end{prop}

\begin{cor}
The mapping $g|M:M\ar \R^k$ is an immersion (i.e. $\rank Dg|M\equiv n-l$)  if and only if at each $p\in M$ 
$$\rank
\left [
\begin{array}{ccc}
   \frac{\partial g_1}{\partial x_1}(p) &\ldots & \frac{\partial g_1}{\partial x_n}(p)\\
   \vdots & \ddots & \vdots\\
   \frac{\partial g_k}{\partial x_1}(p) &\ldots & \frac{\partial g_k}{\partial x_n}(p)\\
   \frac{\partial h_1}{\partial x_1}(p) &\ldots & \frac{\partial h_1}{\partial x_n}(p)\\
   \vdots & \ddots & \vdots\\
   \frac{\partial h_l}{\partial x_1}(p) &\ldots & \frac{\partial h_l}{\partial x_n}(p)\\
\end{array} \right ]
=n,$$
i.e. this matrix has a non-zero $(n\times n)$--minor.
\end{cor}

\begin{cor} \label{immersion}
Let $g=(g_1,g_2,\ldots ,g_{2n}):\R^{n+1},0\ar \R^{2n},0$ be a $C^1$ function. Let us denote $\omega(x):=x_1^2+x_2^2+\ldots +x_{n+1}^2$ for $x\in \R^{n+1}$. Of course, $S^n(r)=\{ x\ | \ \omega (x)-r^2=0\} $.

Then the following conditions are equivalent:
\begin{itemize}
\item[(a)] there exists $r_0>0$ such that $g$ restricted to each sphere of a radius $0<r\leqslant r_0$ is an immersion;
\item[(b)] if  $M_1(x), \ldots ,M_N(x)$ are all the $(n+1)\times (n+1)$--minors of the matrix
$$\left [
\begin{array}{cccc}
    \frac{\partial g_1}{\partial x_1} & \frac{\partial g_1}{\partial x_2} & \ldots & \frac{\partial g_1}{\partial x_{n+1}}\\
    \vdots & \vdots & \ddots & \vdots\\
    \frac{\partial g_{2n}}{\partial x_1} & \frac{\partial g_{2n}}{\partial x_2} & \ldots &  \frac{\partial g_{2n}}{\partial x_{n+1}}\\
    \frac{\partial \omega}{\partial x_1} & \frac{\partial \omega}{\partial x_2} & \ldots & \frac{\partial \omega}{\partial x_{n+1}}\\
\end{array} \right ],$$
then there exists $r_0>0$ such that for each $x\in \R ^{n+1}$ with $0<\| x\|<r_0$ there exists $i\in \{ 1,\ldots ,N\}$ such that $M_i(x)\neq 0$. 
\end{itemize}
\end{cor}

\bigskip
Let $g:\R ^{n+1}\ar \R ^{2n}$ be a $C^2$ mapping. Assume that for each $r\in [r_1,r_2]$, where $0<r_1<r_2$, $g|{S^n(r)}$ is an immersion. We can define a regular homotopy
$$H: [0;1]\times S^n(r_1) \ar \R ^{2n},$$
$$(t,x)\mapsto g\left ( \left (1-t+\frac{tr_2}{r_1}\right )x\right )$$
between $H(0,x)=g|{S^n(r_1)}(x)$ and $H(1,x)=g\left (\frac{r_2}{r_1}x\right )$. By Theorem~\ref{smale}, $\I(g|{S^n(r_1)})=\I(g(\frac{r_2}{r_1}x)|{S^n(r_1)})=\I(g|{S^n(r_2)})$.

If there exists $r_0>0$ such that $g$ restricted to a sphere of a radius $0<r\leqslant r_0$ is an immersion then the intersection number of each $g|{S^n(r)}$ is defined and does not depend on $r$.

\bigskip 

Let $n$ be an even positive integer and let
$$g=(g_1,\ldots,g_{2n}) :\R ^{n+1} \ar \R ^{2n}$$
be a $C^2$ mapping. Suppose that $g$ restricted to each sphere of a radius small enough is a completely regular immersion. 

Let $s\in \R$ , $s\neq 0$, and consider a new mapping $g_s=(sg_1, g_2,\ldots ,g_{2n})$. Then $g_s$ restricted to a sphere of a radius small enough is an immersion. 

Indeed, let  $\tilde{M_1}, \ldots ,\tilde{M}_N$ be all the $(n+1)\times (n+1)$--minors of the matrix:
$$\left [
\begin{array}{cccc}
   s\frac{\partial g_1}{\partial x_1} & s\frac{\partial g_1}{\partial x_2} & \ldots &
 s\frac{\partial g_1}{\partial x_{n+1}}\\
 \vdots & \vdots  & \ddots & \vdots\\
 \frac{\partial g_{2n}}{\partial x_1} & \frac{\partial g_{2n}}{\partial x_2} & \ldots & \frac{\partial g_{2n}}{\partial x_{n+1}}\\
\frac{\partial \omega}{\partial x_1} & \frac{\partial \omega}{\partial x_2} & \ldots & \frac{\partial \omega}{\partial x_{n+1}}\\
\end{array} \right ].$$
Then either $\tilde{M_j} = sM_j$ or $\tilde{M_j} = M_j$, where $M_j$ is the appropriate minor corresponding to the function $g$. By Corollary~\ref{immersion}, there exists $r_0>0$ such that for each $0<\| x\|<r_0$ there exists such $i$ that $M_i(x)\neq 0$, and so $\tilde{M_i}(x)\neq 0$.
It means that for each $s\in \R\setminus \{ 0 \}$ the mapping $g_s$ restricted to a sphere of a radius small enough is an immersion.

Let $r>0$ be small and let $p,q \in S^n(r)$. Then $g(p)=g(q)$ if and only if $g_s(p)=g_s(q)$.
Since $g$ restricted to each sphere of a small radius is completely regular, so is $g_s$. 

\begin{cor} \label{family} $\I (g_s|{S^n(r)})=\sgn (s)\I (g|{S^n(r)})$.
\end{cor}
{\em Proof.} Take $p,q\in S^n(r)$ such that $p\neq q$ and $g(p)=g(q)$. Suppose that $v_1,\ldots,v_n$ form a well--oriented basis in $T_pS^n(r)$, and $w_1,\ldots,w_n$ form a well--oriented basis in $T_qS^n(r)$. 

For any vector $w$, the first coordinate of $Dg_s(p)w$ is equal to the first coordinate of $Dg(p)w$ multiplied by $s$, the other coordinates of $Dg_s(p)w$ are the same as the appropriate coordinates of $Dg(p)w$. Consider two matrices: $$[Dg_s(p)v_1\ldots,Dg_s(p)v_n,Dg_s(q)w_1,\ldots,Dg_s(q)w_n]$$ and $$[Dg(p)v_1\ldots,Dg(p)v_n,Dg(q)w_1,\ldots,Dg(q)w_n].$$ The first row in the first matrix is equal to the first row in the second matrix muptiplied by $s$, other respective rows are identical. Then 
\begin{align*}
\det [Dg_s(p)v_1\ldots,Dg_s(p)v_n,Dg_s(q)w_1,\ldots,Dg_s(q)w_n] \\ 
=s\cdot \det [Dg(p)v_1\ldots,Dg(p)v_n,Dg(q)w_1,\ldots,Dg(q)w_n].
\end{align*}
Hence 
\begin{align*}
& \I(g_s|{S^n(r)})  \\ 
& =\frac{1}{2}\sum \sgn (\det [Dg_s(p)v_1\ldots,Dg_s(p)v_n,Dg_s(q)w_1,\ldots,Dg_s(q)w_n]) \\ 
& =\sgn (s) \frac{1}{2} \sum \sgn (\det [Dg(p)v_1\ldots,Dg(p)v_n,Dg(q)w_1,\ldots,Dg(q)w_n])\\ 
& =\sgn (s)I(g|{S^n(r)}),
\end{align*} 
where $p,q\in S^n(r)$, such that $p\neq q$ and $g_s(p)=g_s(q)$.

\hfill $\Box$

Hence if $g$ restricted to a sphere of a radius small enough is a completely regular immersion, then we have a family $\{ g_s\}$ of $C^2$ mappings and the intersection number of their restrictions to a sphere of a small radius satisfies:
$$\R\setminus\{0\} \ni s \mapsto \I(g_s|{S^n(r)})=\sgn (s) \I(g|{S^n(r)})\in \Z,$$
which is of course determined by a sign of a polynomial. In particular the function $s \mapsto \I(g_s|{S^n(r)})$ is algebraically constructible.

\bigskip

Now let us consider the mapping $$g=(g_1,g_2,g_3,g_4): \R^3 \ar \R^4$$ $$g(x,y,z)=(x,y,xz,yz).$$
Then for each $r>0$, $g|{S^2(r)}$ is an immersion because the matrix:
$$\left [
\begin{array}{ccc}
   \frac{\partial g_1}{\partial x_1} & \frac{\partial g_1}{\partial x_2} & \frac{\partial g_1}{\partial x_3}\\
   \ldots & \ldots & \ldots\\
   \frac{\partial g_4}{\partial x_1} & \frac{\partial g_4}{\partial x_2} & \frac{\partial g_4}{\partial x_3}\\
\frac{\partial \omega}{\partial x_1} & \frac{\partial \omega}{\partial x_2} & \frac{\partial \omega}{\partial x_3}\\
\end{array}  \right ]
=\left [
\begin{array}{ccc}
  1 & 0 & 0\\
  0 & 1 & 0\\
  z & 0 & x\\
  0 & z & y\\
  2x & 2y & 2z\\
\end{array} \right ]$$ has a non-zero $(3\times 3)$--minor at each point $p\in \R^3 \setminus \{0\}$. It is easy to verify that  $g|{S^2(r)}(p)=g|{S^2(r)}(q)$ if and only if $p=(0,0,r)$ and $q=(0,0,-r)$.

At $p=(0,0,r)$ 
$$Dg(p)=
\left [
\begin{array}{ccc}
  1 & 0 & 0\\
  0 & 1 & 0\\
  r & 0 & 0\\
  0 & r & 0\\
\end{array}  \right ]$$
and vectors $v_1=(1,0,0)$ and $v_2=(0,1,0)$ form a well-oriented basis in $T_pS^2(r)$. 

At $q=(0,0,-r)$ 
$$Dg(q)=
\left [
\begin{array}{ccc}
  1 & 0 & 0\\
  0 & 1 & 0\\
  -r & 0 & 0\\
  0 & -r & 0\\
\end{array}  \right ]$$ 
and vectors $w_1=(1,0,0)$ and $w_2=(0,-1,0)$ form a well-oriented basis in $T_qS^2(r)$.

Since
$$\det [Dg(p)v_1,Dg(p)v_2,Dg(q)w_1,Dg(q)w_2]=\det \left [
\begin{array}{cccc}
  1 & 0 & 1 & 0\\
  0 & 1 & 0 & -1\\
  r & 0 & -r & 0\\
  0 & r & 0 & r\\
\end{array}  \right ]=-4r^2,$$
$g|{S^2(r)}$ is a completely regular immersion and $$\I(g|{S^2(r)})=-1.$$
Then for each $s\in \R\setminus \{ 0\}$, the mapping $g_{s}=(sx,y,xz,yz)$ restricted to $S^2(r)$ is a completely regular immersion, and 
$$\I(g_{s}|{S^2(r)})=-\sgn (s).$$
So $\R\setminus \{ 0\} \ni s \mapsto \I(g_{s}|{S^2(r)})=-\sgn s\in \{ -1,1\}$ is an algebraically constuctible function, which is of course nontrivial and constant modulo $2$, but not constant modulo $4$.

\section{Families of analytic mappings} \label{sec:5}

Let $\Lambda \subset \R^p$ be an analytic set.
Let
$$g=(g_1,\ldots ,g_{2n}):\R ^{n+1}\times \Lambda \ar \R ^{2n},$$
where $g_1,\ldots ,g_{2n}$ are analytic functions.

For fixed $\lambda \in \Lambda$ we will denote by $g_{\lambda}: \R ^{n+1} \ar \R ^{2n}$ the mapping defined by $g_{\lambda}(x)=g(x,\lambda)$.

\medskip

\noindent
Let $g_{2n+1}(x)=\omega(x):=x_1^2+x_2^2+\ldots +x_{n+1}^2$, and let
\begin{align*}
G_i(x,y,\lambda) & :=g_i(x,\lambda)-g_i(y,\lambda), \qquad (1\leqslant i \leqslant 2n) \\
G_{2n+1}(x,y) & :=g_{2n+1}(x)-g_{2n+1}(y)=\|x\| ^2-\| y\| ^2, 
\end{align*}
$$G=(G_1,\ldots ,G_{2n}):\R ^{n+1}\times \R ^{n+1}\times \Lambda \ar \R^{2n}.$$
Then $G_1,\ldots ,G_{2n+1}$ are analytic, and there
exist analytic functions $h_{ij}$ such that
$$G_i(x,y,\lambda)=h_{i1}(x,y,\lambda)(x_1-y_1)+
\ldots +h_{i(n+1)}(x,y,\lambda)(x_{n+1}-y_{n+1}).$$ The functions $h_{ij}$ are not uniquely determined. 

We fix such functions $h_{ij}$, and for $1\leqslant i_1<i_2<\ldots <i_{n+1} \leqslant 2n+1$ we define
$$W_{i_1\ldots i_{n+1}}=
\left | \begin{array}{cccc}
h_{i_11} & h_{i_12} & \ldots & h_{i_1(n+1)}\\
h_{i_21} & h_{i_22} & \ldots & h_{i_2(n+1)}\\
\vdots & \vdots & \vdots & \vdots\\
h_{i_{n+1}1} & h_{i_{n+1}2} & \ldots & h_{i_{n+1}(n+1)}
\end{array} \right | $$
By Cramer's rule
$$(x_1-y_1)W_{i_1\ldots i_{n+1}}=
\left | \begin{array}{cccc}
G_{i_1} & h_{i_12} & \ldots & h_{i_1(n+1)}\\
G_{i_2} & h_{i_22} & \ldots & h_{i_2(n+1)}\\
\vdots & \vdots & \vdots & \vdots\\
G_{i_{n+1}} & h_{i_{n+1}2} & \ldots & h_{i_{n+1}(n+1)}
\end{array} \right | $$

$$ \vdots $$

$$(x_{n+1}-y_{n+1})W_{i_1\ldots i_{n+1}}=
\left | \begin{array}{cccc}
h_{i_11} & h_{i_12} & \ldots & G_{i_1}\\
h_{i_21} & h_{i_22} & \ldots & G_{i_2}\\
\vdots & \vdots & \vdots & \vdots\\
h_{i_{n+1}1} & h_{i_{n+1}2} & \ldots & G_{i_{n+1}}
\end{array} \right | $$

\medskip

Points $x$ and $y$ lie on the same sphere in $\R^{n+1}$ centered
at the origin if and only if $\|x\|^2=\|y\|^2$, i.e. if
$G_{2n+1}(x,y)=0$. If that is the case, then $g(x,\lambda)=g(y,\lambda)$ if and only
if $G_1(x,y,\lambda)=\ldots =G_{2n}(x,y,\lambda)=0$.

\medskip
Let us define
\begin{align*}
A & :=\{ (x,y,\lambda)\ |\ G_1(x,y,\lambda)=\ldots =G_{2n}(x,y,\lambda)=G_{2n+1}(x,y)=0\} \\
  & =\{ (x,y,\lambda)\ |\ \exists _{r>0}\ x,y\in S^n(r),\ g(x,\lambda)=g(y,\lambda)\}\cup \{ 0\}\times\{ 0\}\times \Lambda, \\
A_{\lambda} & :=\{ (x,y)\ |\ (x,y,\lambda)\in A\}.
\end{align*}
Then $A$ is closed in $\R^{n+1}\times \R^{n+1}\times \Lambda$.

\medskip

Let $\Delta :=\{ (x,x)\ |\ x\in \R^{n+1} \}$.
Then $\Delta \subset A_{\lambda}$ and $\Delta \times \Lambda$ is closed in $\R^{n+1} \times \R^{n+1} \times \Lambda$. Moreover $(x,y,\lambda)\in \Delta \times \Lambda$ if and only if $x_1-y_1=x_2-y_2=\ldots =x_{n+1}-y_{n+1}=0$.

If $(x,y)\in A_{\lambda}\setminus \Delta$ then $W_{i_1\ldots i_{n+1}}(x,y,\lambda)=0$, for every $1\leqslant i_1<\ldots <i_{n+1} \leqslant 2n+1$.

\medskip
Let us define
\begin{align*}
B & :=A\cap \bigcap  \{ (x,y,\lambda)\ |\ W_{i_1\ldots i_{n+1}}(x,y,\lambda)=0\}, \\
B_{\lambda} & :=A_{\lambda}\cap \bigcap \{ (x,y)\ |\ W_{i_1\ldots i_{n+1}}(x,y,\lambda)=0\} \\ 
& =\{ (x,y)\ |\ (x,y,\lambda)\in B\}.
\end{align*}
Then $B$ is closed in $\R^{n+1}\times \R^{n+1}\times \Lambda$. 
If $(x,y)\in A_{\lambda}$ and $x\neq y$, then $(x,y)\in \bigcap \{ (x,y)\ |\
W_{i_1\ldots i_{n+1}}(x,y,\lambda)=0\}$, so
\begin{equation} \label{ab}
A_{\lambda}\setminus \Delta = B_{\lambda}\setminus \Delta.
\end{equation}

\medskip

For $1\leqslant i,r \leqslant n+1$, if
$$\frac{\partial}{\partial z_i}:=\frac{1}{2}\left ( \frac{\partial}{\partial x_i}-
\frac{\partial}{\partial y_i}\right ),$$
then
$$\frac{\partial}{\partial z_i}(x_r-y_r)=\begin{cases} 0, & i\neq r\\ 1, & i=r \end{cases}.$$
For $1\leqslant j \leqslant 2n+1$,
$$\frac{\partial G_j}{\partial z_i}=\left ( \sum _{r=1}^{n+1}\frac{\partial h_{jr}}
{\partial z_i}(x_r-y_r)\right )+h_{ji};$$
$$\frac{\partial G_j}{\partial z_i}(x,x,\lambda)=h_{ji}(x,x,\lambda).$$
On the other hand we have
\begin{align*}
\frac{\partial G_j}{\partial z_i}(x,y,\lambda) & =\frac{1}{2}\left ( \frac{\partial g_j}{\partial x_i}(x,
\lambda)+\frac{\partial g_j}{\partial x_i}(y,\lambda)\right ); \\
\frac{\partial G_j}{\partial z_i}(x,x,\lambda) & =\frac{\partial g_j}{\partial x_i}(x,\lambda).
\end{align*}
Thus
$$h_{ji}(x,x,\lambda)=\frac{\partial g_j}{\partial x_i}(x,\lambda).$$

\begin{lemma} \label{flambda_imm} Let $\lambda \in \Lambda$. Then  $B_{\lambda}\cap \Delta=\{ 0\}$ in some neighbourhood of the origin  if and only if there exists $r(\lambda)>0$ such that for all $0<r<r(\lambda)$ the mapping $g_{\lambda}=(g_1(\cdot,\lambda),\ldots ,g_{2n}(\cdot,\lambda))$ restricted to a sphere of the radius $r$ is an immersion.
\end{lemma}

\noindent {\it Proof. } We have $h_{ji}(x,x,\lambda)=\frac{\partial g_j}{\partial x_i}(x,\lambda)$, so 
for $(x,x)\in \Delta$ the determinants $W_{i_1\ldots i_{n+1}}(x,x)$ are the $(n+1)\times(n+1)$--minors of the matrix
$$\left [
\begin{array}{cccc}
    \frac{\partial g_1}{\partial x_1} & \frac{\partial g_1}{\partial x_2} & \ldots & \frac{\partial g_1}{\partial x_{n+1}}\\
    \vdots & \vdots & \ddots & \vdots\\
    \frac{\partial g_{2n}}{\partial x_1} & \frac{\partial g_{2n}}{\partial x_2} & \ldots &  \frac{\partial g_{2n}}{\partial x_{n+1}}\\
    \frac{\partial \omega}{\partial x_1} & \frac{\partial \omega}{\partial x_2} & \ldots & \frac{\partial \omega}{\partial x_{n+1}}\\
\end{array} \right ].$$
The function $g_{\lambda}$ satisfies the condition {\em (b)} of Corollary~\ref{immersion} if and only if there exists $r(\lambda)>0$ such that for $0<\|x\|<r(\lambda)$ we have $(x,x)\not \in B_{\lambda}$, i.e. $B_{\lambda}\cap
\Delta=\{ 0\}$ in some neighbourhood of the origin. 

\hfill $\Box$

\medskip

Let us define
$$\delta:\R ^{n+1}\times \R ^{n+1} \ar \R,$$
$$\delta(x,y)=\| x-y\| ^2.$$
Then $\delta \geqslant 0$ and $\delta \inv (0)=\Delta$.

\begin{prop} \label{zero-sets}
If $\lambda \in \Lambda$ then the following conditions are
equivalent:
\begin{itemize}
\item[(a)] there exists $r(\lambda )>0$ such that for all $0<r<r(\lambda)$ the mapping $g_{\lambda}$  restricted to a sphere of the radius $r$ is an immersion;
\item[(b)] $\{ 0\}$ is isolated in the set $$\Delta \cap \bigcap W_{i_1\ldots
i_{n+1}}(\cdot, \lambda)\inv (0)=\Delta \cap \bigcap W_{i_1\ldots
i_{n+1}}(\cdot, \lambda)\inv (0)\cap A_{\lambda},$$ where $1\leqslant i_1<\ldots <i_{n+1}\leqslant 2n+1$;

\item[(c)] $B_{\lambda}\cap \Delta=\{ 0\}$ in a neighbourhood of the origin.
\end{itemize}
If that is the case, then 

\begin{itemize}
\item[(d)] in a neighbourhood of the origin
$\overline{A_{\lambda}\setminus \Delta}=B_{\lambda}$, and $(A_{\lambda}\setminus \Delta) \cup \{ 0\}$ is closed;
\item[(e)] if $r>0$ is small enough then $S^{n}({r}/{\sqrt{2}})\times S^{n}({r}/{\sqrt{2}})\cap (G(\cdot,\cdot,\lambda))\inv (0)\setminus \Delta$ is a compact subset of $S^{n}({r}/{\sqrt{2}})\times S^{n}({r}/{\sqrt{2}})\setminus \Delta$.
\end{itemize}
\end{prop}

\noindent {\it Proof. }
Lemma~\ref{flambda_imm} implies that $(a)\Leftrightarrow (c)$. Since $\Delta \subset A_{\lambda}$,  $$\Delta \cap \bigcap
_{1\leqslant i_1<\ldots <i_{n+1} \leqslant 2n+1}W_{i_1\ldots i_{n+1}}(\cdot, \lambda)\inv (0)\cap A_{\lambda}=\Delta \cap B_{\lambda},$$  so $(b)\Leftrightarrow (c)$. By (\ref{ab}) $A_{\lambda}\setminus \Delta = B_{\lambda}\setminus \Delta$,  so $\overline{A_{\lambda}\setminus \Delta} = \overline{B_{\lambda}\setminus \Delta}$. If $B_{\lambda}\cap \Delta=\{ 0\}$ in some neighbourhood of the origin, then $\overline{B_{\lambda}\setminus \Delta}=B_{\lambda}$ in some neighbourhood of the origin, hence $(c)\Rightarrow (d)$. 

By (\ref{ab}), $A_{\lambda}\setminus \Delta=B_{\lambda}\setminus \Delta$ and by $(c)$ (maybe after making $r(\lambda)$ smaller) for $0<r<r(\lambda )$ we have $B_{\lambda }\cap S^{2n+1}(r)\cap \Delta=\emptyset$, so
\begin{align*}
B_{\lambda }\cap S^{2n+1}(r) & =(B_{\lambda }\setminus \Delta )\cap S^{2n+1}(r)=(A_{\lambda }\setminus \Delta )\cap S^{2n+1}(r) \\ 
& =\{ (x,y)\in S^{2n+1}(r) \ |\ \|x\|^2=\|y\|^2,\ g_{\lambda}(x)=g_{\lambda }(y),\ x\neq y \} \\
& =(S^{2n+1}(r)\cap G_{2n+1}\inv (0)\cap (G(\cdot, \cdot, \lambda))\inv (0))\setminus \Delta \\
& =S^{n}({r}/{\sqrt{2}})\times S^{n}({r}/{\sqrt{2}})\cap (G(\cdot, \cdot, \lambda))\inv (0))\setminus \Delta.
\end{align*}
$B_{\lambda}$ is closed, so $B_{\lambda }\cap S^{2n+1}(r)$ is compact and we have $(e)$.

\hfill $\Box$

\section{Immersions of spheres --- the algebraic case} \label{algebraic} \label{sec:6}

Assume that $\Lambda \subset \R^p$ is an algebraic set and each $g_i$ is a polynomial. Then there
exist polynomials $h_{ij}$ such that
$$G_i(x,y,\lambda)=h_{i1}(x,y,\lambda)(x_1-y_1)+
\ldots +h_{i(n+1)}(x,y,\lambda)(x_{n+1}-y_{n+1}).$$ Now $W_{i_1\ldots i_{n+1}}$ are polynomials, and the sets $A$, $B$, $A_{\lambda}$ and $B_{\lambda}$ are algebraic.

\medskip
Let us define
\begin{align*}
Z & =(\delta \inv (0)\times \Lambda )\cap B=(\Delta\times \Lambda )\cap B, \\
Z_{\lambda} & =\delta \inv (0) \cap B_{\lambda}=\Delta \cap B_{\lambda}.
\end{align*}
So $\delta (x,y,\lambda):=\delta(x,y)$ is continuous on $B$
and the assumptions of \cite[Corollary 3.1]{fekak} hold. Then there exist a finite division $\Lambda =\bigcup S_i$ (where $S_i$ are semialgebraic), continuous, semialgebraic functions $h_i:B\cap(\R ^{2n+2}\times {S_i})\ar \R$, and constants ${q_i}\in \Q ^+$
such that for
$\lambda \in S_i$, $(x,y)\in B_{\lambda}$
$$\dist ((x,y), Z_{\lambda})^{{q_i}}\leqslant h_i(x,y,\lambda)\delta(x,y).$$

If $\lambda \in S_i$, then there exists such a constant $c_{\lambda}>0$ that  $$|h_i(x,y,\lambda)|<c_{\lambda}$$
for each $(x,y)\in B_{\lambda}$ sufficiently close to the origin.

\medskip

Set $d(x,y)=x_1^2+\ldots +x_{n+1}^2+y_1^2+\ldots +y_{n+1}^2=\| (x,y)\| ^2$. From now on we assume that for each $\lambda \in \Lambda$ there exists $r(\lambda)>0$ such that for each $0<r<r(\lambda)$ the mapping $g_{\lambda}|S^n(r)$ is an immersion. 

By Proposition~\ref{zero-sets}, locally $Z_{\lambda}=\Delta\cap B_{\lambda}=\{ 0\}$. Then for $(x,y)\in B_{\lambda}$ close to the origin
$$d(x,y)^{\frac{q_i}{2}}=\dist ((x,y),\{ 0\})^{{q_i}}=\dist ((x,y),Z_{\lambda})^{{q_i}}\leqslant c_{\lambda}\delta(x,y).$$
If $\alpha > \max _i\{ \frac{q_i}{2}\} $ is an integer, then for each $\lambda \in \Lambda$ and $t>0$
\begin{equation} \label{LI} \delta (x,y)\geqslant td(x,y)^{\alpha}\end{equation}
for $(x,y)\in B_{\lambda}$ lying sufficiently close to the origin. By (\ref{ab}) $A_{\lambda}\setminus \Delta =B_{\lambda}\setminus \Delta$, so the inequality (\ref{LI}) holds on $A_{\lambda}\setminus \Delta$, i.e. on 
$$\{ (x,y)\ |\ G_1(x,y,\lambda)=\ldots =G_{2n}(x,y,\lambda)=\|x\|^2-\|y\|^2=0\} \setminus \delta \inv (0).$$
By Proposition~\ref{zero-sets} (d), $(A_{\lambda}\setminus \Delta) \cup \{0\}$ is closed.  

\medskip

A polynomial $F_1(x,y)=\| x\| ^2-\| y\| ^2$  has an isolated critical point at the origin. 
Then
$$L(r)=S^{2n-1}(r)\cap F_1\inv (0)=S^{n}({r}/{\sqrt{2}})\times S^{n}({r}/{\sqrt{2}})$$
is a compact oriented $2n$--dimensional manifold. Since $\delta (x,y)=\|x-y\| ^2\geqslant 0$,
$$U(r)=L(r)\setminus \delta \inv (0)=S^{n}({r}/{\sqrt{2}})\times S^{n}({r}/{\sqrt{2}})\setminus \Delta.$$

By Proposition~\ref{isolsing} there exists $\alpha _0>0$ such that for any integer $\alpha > \alpha _0$ and $t\neq 0$, $(F_1, \delta -td^{\alpha})\inv (0)$ has an isolated singularity at the origin. (As $F_1$, $\delta$ and $d$ are homogeneous of degree $2$, it is easy to verify that $\alpha _0=1$ would be enough.)

\medskip

Let us assume that $n$ is even.
Take an integer $\alpha > \max (\max _i\{ \frac{q_i}{2}\} ,\alpha _0)$. Put $F_2=\delta - td^{\alpha}$ for arbitrary $t>0$, and 
\begin{align*}
F & =(F_1,F_2), \\
G_{\lambda} & =G(\cdot,\cdot,\lambda):\R^{n+1}\times \R ^{n+1},0\ar \R^{2n},0.
\end{align*}
By Proposition~\ref{zero-sets} (e), $L(r)\cap G_{\lambda}\inv (0)\setminus \delta \inv (0)$ is a compact subset of $U(r)$. By Proposition~\ref{topdeg}, for $r>0$ small enough, the topological degree of
$${\textstyle (S^{n}(\frac{r}{\sqrt{2}})\times S^{n}(\frac{r}{\sqrt{2}})\setminus \Delta, S^{n}(\frac{r}{\sqrt{2}})\times S^{n}(\frac{r}{\sqrt{2}})\setminus (G_{\lambda}\inv (0)\cup \Delta)) \ar (\R ^{2n}, \R ^{2n}\setminus \{ 0\})}$$
$$(x,y) \mapsto G_{\lambda}(x,y)=g_{\lambda}(x)-g_{\lambda}(y),$$
i.e. $2\I (g_{\lambda}|S^n(r/\sqrt{2}))$, equals the local topological degree $\deg _0(F,G_{\lambda})$.

\bigskip

Let us define a polynomial mapping $H:\R \times \R^{n+1}\times \R^{n+1}\times \Lambda \ar \R\times \R\times \R^{2n}$:
$$H(t,x,y,\lambda)=(F_1(x,y),\delta (x,y) -td(x,y)^{\alpha},g(x,\lambda)-g(y,\lambda)).$$
Denote $H^t_{\lambda}(x,y)=H(t,x,y,\lambda)$. For any $\lambda \in \Lambda$ and $t>0$, the mapping $H^t_{\lambda}$ has an isolated zero at the origin and the local topological degree $\deg _0 (H^t_{\lambda})=2\I(g_{\lambda}|{S^n(r)})$ for $0<r<r(\lambda)$. If $t<0$ then $\delta-td^{\alpha}>0$, except of the origin, so that $\deg _0 (H^t_{\lambda})=0$.

\begin{theorem} \label{alg_constr}
Let $\Lambda \subset \R ^p$ be an algebraic set and let $n$ be an even integer.
Let
$$g_1(x,\lambda),\ldots ,g_{2n}(x,\lambda):\R ^{n+1}\times \Lambda \ar \R ^{2n}$$
be polynomials.
Assume that for each $\lambda \in \Lambda$ there exists $r(\lambda)>0$ such that for each $0<r<r(\lambda)$ the mapping $g_{\lambda}$ restricted to a sphere $S^n(r)$ of the radius $r$ centered at the origin is an immersion.
Then the function $\Lambda \ni \lambda \mapsto \I(g_{\lambda}|{S^n(r)})\in \Z$ is algebraically constructible.
\end{theorem}

\noindent {\it Proof. } 
Using the same arguments as in \cite[Theorem 3.4]{szafr-parus1} one may prove that there exist polynomials $h_1,\ldots ,h_s$ on  $\R\times \Lambda$, such that if $0$ is isolated in $(H^t_{\lambda})\inv (0)$ then
$$\deg _0 (H^t_{\lambda})=\sum _{i=1}^{s} \sgn h_i(t,{\lambda}),$$
i.e. this function is algebraically constructible.

For $r>0$ small enough
\begin{align*}
\I(g_{\lambda}|{S^n(r)}) & =\lim _{t\to 0}\frac{1}{2}(\deg _0 (H^t_{\lambda})+\deg _0 (H^{-t}_{\lambda})) \\ 
& =\lim _{t\to 0}\frac{1}{2}\sum _{i=1}^{s} \left (\sgn h _i({\lambda},t)+\sgn h _i({\lambda},-t)\right )=:\psi (\lambda).
\end{align*}
According to \cite[Lemma 6.5]{szafr-parus1} the function $\Lambda \ni \lambda \mapsto \psi(\lambda)\in \Z$ is algebraically constructible, and so is $\Lambda \ni \lambda \mapsto \I(g_{\lambda}|{S^n(r)})$.

\hfill $\Box$

\noindent
{\em Remark.} It is easy to see that if mappings $g_{\lambda}$ are immersions on small spheres only for $\lambda \in \Lambda \setminus \Sigma$, where $\Sigma$ is a proper algebraic subset of $\Lambda$, then the function $\Lambda \ni \lambda \mapsto \I(g_{\lambda}|{S^n(r)})\in \Z$ is generically algebraically constructible, i.e. it coincides with an algebraically constructible function in $\Lambda \setminus \Sigma$.

\end{document}